\documentstyle{amsppt}
\hsize =35pc
\vsize =53pc
\NoBlackBoxes
\topmatter
\title A non-quasiconvex subgroup of a hyperbolic group with an exotic limit set
\endtitle
\rightheadtext {A non-quasiconvex subgroup}
\author Ilya Kapovich \endauthor
\thanks{This research is supported by Alfred P. Sloan Doctoral Dissertation Fellowship}
\endthanks
\address{City College, 138th Street and Convent Avenue, 
New York, NY 10031}
\endaddress
\email ilya\@groups.sci.ccny.cuny.edu\endemail
\abstract We construct an example of a torsion free freely indecomposable finitely presented non-quasiconvex subgroup $H$ of a word hyperbolic group $G$ such that the limit set of $G$ is not the limit set of a quasiconvex subgroup of $G$. In particular, this gives a counterexample to the conjecture of G.Swarup that a finitely presented one-ended subgroup of a word hyperbolic group is quasiconvex if and only if it has finite index in its virtual normalizer.
\endabstract

\endtopmatter

\document

\head 1. Introduction \endhead
 
A subgroup $H$ of a word hyperbolic group $G$ is {\it quasiconvex (or rational)} in $G$ if for any finite generating set $A$ of $G$ there is $\epsilon >0$ such that every geodesic in the Cayley graph $\Gamma(G,A)$ of $G$ with both endpoints in $H$ is contained in $\epsilon$-neighborhood of $H$. The notion of a quasiconvex subgroup corresponds, roughly speaking, to that of geometric finiteness in the theory of classical hyperbolic groups (see [Swa], [KS], [Pi]).
Quasiconvex subgroups of word hyperbolic groups are finitely presentable and word hyperbolic and their finite intersections are again quasiconvex.  Non-quasiconvex finitely generated subgroups of word hyperbolic groups are quite rare and there are very few examples of them. We know only three basic examples of this sort.
The first is based on a remarkable construction of E.Rips [Ri], which allows one, given an arbitrary finitely presented group $Q$, to construct a word hyperbolic group $G$ and a two-generator subgroup $H$ of $G$ such that $H$ is normal in $G$ and the quotient is isomorphic to $Q$. The second example is based on the existence of a closed hyperbolic 3-manifold fibering over a circle, provided by results of W.Thurston and T.Jorgensen. The third example is obtained using the result of M.Bestvina and M.Feign [BF] who proved that if $F$ is a non-abelian free group of finite rank and $\phi$ is an automorphism of $F$ without periodic conjugacy classes then the HNN-extension of $F$ along $\phi$ is word hyperbolic.   

If $G$ is a word hyperbolic group then we denote the {\it boundary} of $G$ (see [Gr], [GH] and [CDP]) by $\partial G$.
For a subgroup $H$ of $G$ the {\it limit set} $\partial_G(H)$ of $H$ is the set of all limits in $\partial G$ of sequences of elements of $H$. 

In this note we construct an example of a non-quasiconvex finitely presented one-ended subgroup $H$ of a word hyperbolic group $G$ such that the limit set of $H$ is {\it exotic}. By exotic we mean that the limit set of $H$ is not the limit set of a quasiconvex subgroup of $G$.
This result is of some interest since in the previously known examples non-quasiconvex subgroups were normal in the ambient hyperbolic groups and thus (see [KS]) had the same limit sets.
Our subgroup $H$ also provides a counter-example to the conjecture of G.A. Swarup [Swa] which stated that a finitely presented freely indecomposable subgroup of a torsion-free word hyperbolic group is quasiconvex if and only if it has finite index in its virtual normalizer (this statement was known to be true for 3-dimensional Kleinian groups). 
The subgroup $H$, constructed here, coincides with its virtual normalizer.
(By the virtual normalizer of a subgroup $H$ of a group $G$ we mean the subgroup $VN_G(H)=\{g\in G\mid |H: H\cap gHg^{-1}|<\infty, |gHg^{-1}:H\cap gHg^{-1}|<\infty\}$.)

\head 2. The proofs \endhead
\proclaim {Proposition 1} 
Let $F$ be the fundamental group of a closed hyperbolic surface $S$.
Let $\phi$ be an automorphism of $F$ induced by a pseudo-anosov homeomorphism of $S$.
Take $G$ to be the mapping-torus group of $\phi$
that is $G=<F,t | tft^{-1}=\phi(f), f\in F>$.
Let $x\in F$ be an element which is not a proper power in $F$
(and so, obviously, $x$ is not a proper power in $G$).
Let $G_1$ be a copy of $G$.
The group $G_1$ contains a copy $F_1$ of $F$ and a copy $x_1$ of $x$.
Put $$M=G \underset {x=x_1}\to\ast G_1 \eqno (1)$$ and $H=sgp(F,F_1)$.

Then 
\roster
\item $M$ is torsion-free and word hyperbolic.
\item $H$ is finitely presented, freely indecomposable
and non-quasiconvex in $M$.
\endroster
\endproclaim

\demo{Proof}
The group $G$ is torsion free and word hyperbolic since
it is the fundamental group of a closed 3-manifold of constant
negative curvature (see [Th]).
Thus $M$ is word hyperbolic by the combination theorem for
negatively curved groups (see [BF], [KM]).
Notice that $H=sgp(F,F_1)\cong F \underset {x=x_1}\to\ast F_1$ and so $H$ is torsion-free,
finitely presentable and freely indecomposable. Moreover, $H$ is word hyperbolic
by the same combination theorem.

Suppose $H$ is quasiconvex in $M$. It is not hard to show (see [Ka])
that $F$ is quasiconvex in $H$ since $H=F\ast_C F_1$ where $C=<x>=<x_1>$ is cyclic.
Thus $F$ is quasiconvex in $M$.
On the other hand $G$ is also quasiconvex in $M$ since $M=G\ast_C G_1$ and $C$ is cyclic.
This implies that $F$ is quasiconvex in $G$ which contradicts the fact
that $F$ is infinite, normal in $G$ and has infinite index in $G$ (see [ABC]).
Thus $H$ is not quasiconvex in $M$.
\enddemo

\proclaim {Theorem 2} Let $G$, $G_1$, $M$ and $H$ be as in Proposition 1.
Put $K$ to be the limit set $\partial_M(H)$ of $H$ in the boundary $\partial(M)$ of $M$.
Then $H=Stab_M(K)=\{f\in M| fK=K\}$.
\endproclaim
\demo{Proof}
Before proceeding with the proof, we choose a finite generating set $\Cal G$ for $G$ and its copy $\Cal G_1$.
Then $\Cal G$ defines the word length $l_G$ on $G$ and the word metric $d_G$ on $G$.
Analogously, $\Cal G \cup \Cal G_1$ is a finite generating set for $M=G\ast_C G_1$ which defines the word length $l_M$ and the word metric $d_M$ on $M$. Fix a $\delta >0$ such that all $d_M$-geodesic triangles are $\delta${\it -thin} (see [ABC]).
We also denote by $C$ the subgroup of $M$ generated by $x=x_1$.
Thus $M=G\underset C\to \ast G_1$ and $H=F\underset C\to\ast F_1$.

Suppose $z\in Stab_M(K)$.
We will show that $z\in H$ by induction on the syllable length of $z$
with respect to presentation (1).
When the syllable length of $z$ is 0, that is $z\in C$, the statement is obvious.
Suppose now that $z\in Stab_M(K)-H$, the syllable length of $z$ is $m>0$ and the statement has
been proved for elements of smaller syllable length.
Write $z$ as a strictly alternating product $z=u_1..u_m$ of elements
from $G-C$ and $G_1-C$.
If $u_m\in F\cup F_1$, then $u_m\in H\cap Stab_M(K)$ and so
$u_1..u_{m-1}\in Stab_M(K)$. Therefore $u_1..u_{m-1}\in H$ by the inductive hypothesis, $u_m\in H$ and so $z\in H$.
Thus $u_m\in (G-F)\cup(G_1-F_1)$.
Assume for definiteness that $u_m\in G_1-F_1$ that is $u_m=f_1t_1^j$ for some $j\not=0$, $f_1\in F_1$.

Choose $y\in F$ so that no power of $y$ is conjugate in $G$
to a power of $x$.
We may further assume that $y$ is periodically geodesic in $G$, that is,
$l_G(y^n)=|n|l_G(y)$ for each $n$.
Indeed, for some $k$ $y^k$ is conjugate to a periodically geodesic element $y'\in G$.
Since $F$ is normal in $G$, $y'\in F$ and we may replace $y$ by $y'$. Fix a $d_G$-geodesic representative $Y$ of $y$.
\enddemo
\proclaim {Claim 1} There is a constant $K$ depending on $x$ and $y$ such that for any integer $n$
there is an element $c\in C$ of $l_G$-length at most  $K$ 
such that $y^n=cu$ where $u$ is the shortest (with respect to $l_G$) element in the coset class $Cy^n$.
\endproclaim
This follows from quasiconvexity of the subgroups $<x>$ and $<y>$ in $G$ and from the fact
that no power of $x$ is conjugate to a power of $y$ in $G$.

\proclaim {Claim 2} There is a constant $D>0$ depending on $x$ and $y$ such that if $q\in M$ either ends in $G_1-C$ or belongs to $C$ then $l_M(qy^n)\ge |n|/D-D$ for any integer $n$.
\endproclaim
To prove Claim 2 recall that by the theorem of G.Baumslag, S. Gersten, M.Shapiro, and H.Short [BGSS], cyclic amalgamations of hyperbolic groups are automatic. In the proof of this theorem they construct an actual automatic language on a cyclic amalgam of two hyperbolic groups, which, therefore, consists of quasigeodesic words [ECHLPT, Theorem 3.3.4].
We will explain how their procedure works in the case of the group $M=G\ast_C G_1$ (we use the fact that $d_G\mid_C=d_{G_1}\mid_C$).
Suppose we have an element $e\in M$ and we want to find its representative in the automatic language $L$. If $e\in C$, then $e=x^k$ and we take a $d_G$-geodesic representative of $x^k$ to be the representative of $e$ in the automatic language on $M$. Suppose $e\not\in C$.
First, write $e$ as a strictly alternating product of elements from $G-C$ and $G_1-C$,
$e=e_1..e_j$. 
Then express $e_1$ as ${e_1}=\overline {w_1}c^{n_1}$ where $w_1$ is a $d_G$-geodesic ($d_{G_1}$-geodesic) word such that $\overline w_1$ is the shortest in the $d_G$-metric ($d_{G_1}$-metric) in the coset class ${e_1} C$.
Then express $c^{n_1}{e_2}$ as  $c^{n_1}{e_2}=\overline {w_2}c^{n_2}$ where $w_2$ is a $d_G$-geodesic ($d_{G_1}$-geodesic) word such that $\overline w_2$ is the shortest in the $d_G$-metric ($d_{G_1}$-metric) in the coset class $c^{n_1}e_2C$. And so on for $i=1,2,..j-1$. For the very last $w_j$ we put $w_j$ to be a $d_G$-geodesic ($d_{G_1}$-geodesic) representative of $c^{n_{j-1}}e_j$.
As a result we obtain the word $w=w_1..w_j$ such that $\overline w=e$. This word $w$ is the required representative of $e$ in $L$.
\smallskip
Let us return to the proof of Claim 2.
Suppose $q\in M$ and either $q$ ends in $G_1-C$ or $q\in C$.
When we apply the above rewriting process to find the representative of $qy^n$ in $L$, we see that by Claim 1 the very last word $w_j$ will be at most $2K$-shorter than the $d_G$-geodesic word $y^n$ (where $K$ is the constant provided by Claim 1).
Since the subgroup $<y>$ is quasiconvex in $M$ and the language $L$, which our rewriting process yields, consists of $d_M$-quasigeodesics, Claim 2 follows. 
\smallskip
Let $y^+=\underset{n\rightarrow\infty}\to\lim y^n\in \partial M$.
By definition of $K$ we have $y^+\in K$ and therefore $zy^+\in K$.
This means that for any $N>0$ there is an element $h\in H$ and a positive power $y^n$ of $y$
such that $(h,zy^n)_1>N$, the Gromov inner product taken in the $d_M$-metric.
This means that $d_M(h(N), (zy^n)(N))\le \delta$ where
$h(N)$ and $(zy^n)(N)$ are elements of $M$ represented by initial segments of length $N$
of $d_M$-geodesic representatives of $h$ and $zy^n$.
We will obtain a contradiction by showing that $l_M(h(N)^{-1}(zy^n)(N))$
is big when $N$ is big.
More precisely, we will show that 
$$l_M(h(N)^{-1}(zy^n)(N))> K_0\cdot N\ - K_0 \eqno (2)$$
where $h\in H$, $n\in {\Bbb Z}$ and $K_0>0$ is a constant independent of $h$, $N$ and $n$.

If $h\in C$ we put $p_1=h^{-1}$.
If $h\not\in C$, we express $h^{-1}$ as an alternating product $h^{-1}=p_1p_2..p_s$ of elements of $F-C$ and $F_1-C$
where for all $i<s$ the element $p_i$ is shortest in the coset class $p_iC$ in the $G$-metric
when $p_i\in F$ and $p_i$ is shortest in the coset class $p_iC$ in the $G_1$-metric
when $p_i\in F_1$. Let $\hat {p_i}$ be a $d_G$-geodesic ($d_{G_1}$-geodesic) representative of $p_i$. Then $W_1=\hat {p_1}....\hat {p_s}$ is a $\lambda$-quasigeodesic representative
of $h^{-1}$ with respect to the $d_M$-metric where $\lambda >0$ is some constant independent of $h$. (see the proof of Claim 2).
Recall that $z=u_1...u_m$, where $u_m=f_1t_1^j$, $j\not=0$.
Let $\hat u_i$ be a $d_G$-geodesic ($d_{G_1}$-geodesic) representative of $u_i$, $i\le m$.
Then $W_2=\hat u_1\dots\hat u_{m-1}\hat u_m Y^n$ is a $\lambda_1$-quasigeodesic representative of $zy^n$ with respect to $d_M$,
where the constant $\lambda_1$ does not depend on $n$ (note that the element $z$ is fixed).
Put $\lambda_2=max(\lambda, \lambda_1)$. Let $\epsilon>0$ be  such that any two $\lambda_2$-quasigeodesics with common endpoints in the Cayley graph of $M$ are $\epsilon$-Hausdorff-close.
Then for an element $h(N)$ there is a terminal segment $T_1$ of $W_1$ such that $l_M(\overline {T_1}\cdot h(N))\le \epsilon$.
Likewise, for an element $(zy^n)(N)$ there is an initial segment $T_2$ of $W_2$ such that $l_M(\overline {T_2}^{-1}\cdot (zy^n)(N))\le \epsilon$.
\smallskip
Thus to prove (2) it is enough to show that if $S_1$ is a terminal segment of $W_1$ of
length $N_1$ and $S_2$ is an initial segment of $W_2$ of length $N_2$
then $$L_M(\overline {S_1S_2})\ge K_1\cdot min(N_1,N_2)\ -K_1 \eqno (3)$$
where $K_1>0$ is some constant independent og $n$ and $h$.
\smallskip
 Since $z$ and $y$ are fixed, we may assume that $S_2=\hat u_1..\hat u_{m} Y^k$ for some $0\le k\le n$.
Then $S_1$ has the form  $S_1=\hat q_i\hat p_{i+1}..\hat p_s$ where $i\le s$ and $\hat q_i$ is a nonempty terminal segment of $\hat p_i$. Put $q_i=\overline {\hat q_i}$.
\smallskip
\proclaim {Claim 3} Either $q_ip_{i+1}..p_su_1..u_m$ ends (when rewritten in the normal form with respect to (1)) in the element
of $G_1-C$ or $q_ip_{i+1}..p_su_1..u_m\in C$.
\endproclaim

Indeed, $u_m\in G_1-C$ and so $q_ip_{i+1}..p_su_1..u_m$ ends (when rewritten in the normal form with respect to (1)) in the element
of $G-C$ unless $u_m^{-1}..u_1^{-1}$ is a {\it terminal segment} of $q_ip_i..p_s$ that is 
either $m>s-i+1$ and for some $i\le j\le m$
$$p_j..p_su_1..u_m\in C \eqno (4)$$
or 
$$q_ip_i..p_su_1..u_m\in C  \eqno (5)$$

It is clear that (4) is impossible since $p_j..p_s\in H$ and  $u_1..u_m\not\in H$.
If (5) holds, we have $q_ip_{i+1}..p_su_1..u_m\in C$ as required. Thus Claim 3 is established.
\smallskip
We have verified that either $q_ip_{i+1}..p_su_1..u_m\in C$ or $q_ip_{i+1}..p_su_1..u_m$ ends in the element of $G_1-C$.  
Therefore by Claim 2 $l_M(q_ip_{i+1}..p_su_1..u_my^k)=l_M(\overline {S_1S_2})\ge k/D-D$ and (3) follows.
This completes the proof of Theorem 2.

\proclaim {Corollary 3} Let $M$, $G$, $G_1$, $C$ and $H$ be as in Theorem 2. Then 
\roster
\item "(a)" the limit set of $H$ is not the limit set of a quasiconvex subgroup of $M$;
\item "(b)" the virtual normalizer $VN_M(H)$ of $H$ in $M$ is equal to $H$
\endroster
\endproclaim
\demo{Proof}

(a) Suppose there is a quasiconvex subgroup $Q_1$ of $M$ such that
$\partial_M(H)=\partial_M(Q_1)=K$. Clearly $Q_1$ is infinite since $K$ is nonempty.
Put $Q=Stab_M(K)=\{y\in M | yK=K\}$.
Since $Q_1$ is infinite and quasiconvex in $M$ and $Q=Stab_M(\partial_M(Q_1))$, it follows from [KS, Lemma 3.9] that 
$Q$ contains $Q_1$ as a subgroup of finite index and therefore $Q$ is also quasiconvex in $M$.
On the other hand Theorem 2 implies that $H=Stab_M(K)$ and so $H=Q$. This contradicts the fact that $H$ is not quasiconvex in $M$ by Proposition 1.
\smallskip
(b) It is not hard to see, that $A\le VN_B(A)\le Stab_B(\partial_B(A))$ when $A$ is an infinite subgroup of a word hyperbolic group $B$.
Indeed, if $g\in VN_B(A)$ then $A_0=A\cap gAg^{-1}$ has finite index $n$ in $A$. Let $A=A_0\cup A_0c_1\dots \cup  A_0c_{n-1}$ and let $D=max\{ l_A(c_i)\vert i=1,.., n-1\}$. Suppose $p\in\partial_B(A)$. Then there is a sequence $a_m\in A$ such that $p=\underset{m\rightarrow\infty}\to\lim a_m$. For each $m$ there is $b\in B$ with $l_B(b)\le D+l_B(g)$ such that $ga_mb=a_m'\in A_0$. Therefore $gp\in \partial_B(A_0)=\partial_B(A)$. Since $p\in\partial_B(A)$ was chosen arbitrarily, we have $g\partial_B(A)\subseteq\partial_B(A)$. Since by the same argument $g^{-1}\partial_B(A)\subseteq\partial_B(A)$ we conclude that $g\partial_B(A)=\partial_B(A)$. Thus $A\le VN_B(A)\le Stab_B(\partial_B(A))$.

 For the subgroup $H$ of $M$ we have $H\le VN_M(H)\le Stab_M(\partial_M(H))$. On the other hand $Stab_M(\partial_M(H))=H$ by Theorem 2. Therefore $H=VN_M(H)$.

\enddemo

\enddocument